\documentclass[12pt]{article}
\usepackage{latexsym,fullpage}
\usepackage{graphics,graphicx}
\usepackage{subfigure}
\usepackage{psfrag}
\usepackage{epic,eepic,color,verbatim}
\usepackage{bbm,euscript,amssymb,amsthm}
\usepackage{amsfonts}
\usepackage{amsmath}
\usepackage{amssymb}
\usepackage{amscd}
\usepackage{amsthm}
\usepackage{mathrsfs}
\usepackage{cite}
\usepackage{graphicx}
\usepackage{pstricks, pstricks-add, pst-math}
\usepackage{pst-xkey}
\usepackage{color}
\usepackage{calc}

\usepackage{verbatim}

\usepackage{soul}

\usepackage{authblk}

\def\mytextindent#1{\indent\indent\llap{\rm#1\enspace}\ignorespaces}
\def\myitem{\par\hangindent0pt\mytextindent}
\def\myclaim#1#2{\bigskip\noindent\rlap{\rm(#1)}\ignorespaces
   \hangindent=30pt\hskip30pt{\sl#2}\bigskip}
\def\mylabel#1{{\label{#1}}}
\def\proof{\noindent {\bf Proof. }}

\def\ignore#1{{}}

\DeclareSymbolFont{slantedbx}{OT1}{cmr}{bx}{sl}
\DeclareSymbolFontAlphabet\mathslbx{slantedbx}
\DeclareSymbolFont{AMSa}{U}{msa}{m}{n} \DeclareMathSymbol{\upharpoonright} {\mathrel}{AMSa}{"16}

\newtheorem{thm}{Theorem}[section]

\def\kr{{\rm cr}}

\newtheorem{theorem}[thm]{Theorem}

\newtheorem{lemma}[thm]{Lemma}

\newcounter{claim}[thm]

\def\scr#1#2{{\mathop{\rm cr_{#1}}(#2)}}
\def\crg(#1){\scr{}{#1}}
\def\crgp(#1){\scr{P}{#1}}

\def\kr{{\hbox{\rm cr}}}

\def\comm#1{{\marginpar{\tiny #1}}}

\title{\bf Nested Cycles in Large Triangulations and Crossing-Critical Graphs}

\author[1]{C\'esar Hern\'andez--V\'elez}
\author[,1]{Gelasio Salazar\thanks{Partially supported by CONACYT Grant 106432.}}
\author[,2]{Robin Thomas\thanks{{Partially supported by NSF grant DMS-0739366.}}}
\affil[1]{Instituto de F\'\i sica, Universidad Aut\'onoma de San Luis
  Potos\'{\i}. San Luis Potos\'{\i}, Mexico 78000}
\affil[2]{School of Mathematics, Georgia Institute of
  Technology. Atlanta, GA, 30332}

\begin{document}
\maketitle

\begin{abstract}
We show that every sufficiently large plane triangulation has a
large collection of nested cycles that either are pairwise disjoint,
or pairwise intersect in exactly one vertex, or pairwise intersect in
exactly two vertices.  We apply this result to show that for each fixed
positive integer $k$, there are only finitely many $k$-crossing-critical 
simple graphs of average degree at least six. 
Combined with the recent constructions of
crossing-critical graphs given by Bokal, this settles the question of
for which numbers $q>0$ there is an infinite family of
$k$-crossing-critical simple graphs of average degree~$q$.
\end{abstract}

\section{Introduction}

All {\em graphs} in this paper are finite, and may have loops and parallel edges.
The {\em crossing number} of a graph $G$, denoted by $\kr(G)$, is
the minimum, over all drawings $\gamma$ of $G$ in the plane,
of the number of crossings in $\gamma$.
(We will formalize the notion of a drawing later.)

A graph $G$ is $k$-{\em crossing-critical} if the crossing number
of $G$ is {at least} $k$ and $\kr(G-e) < k$ for every edge $e$ of
$G$. The study of crossing-critical graphs is a central part of the
emerging structural theory of crossing numbers. Good examples of this
aspect of crossing numbers are Hlin\v{e}n\'y's proof that
$k$-crossing-critical graphs have bounded path-width~\cite{petr};
Fox and T\'oth's work on the decay of crossing numbers~\cite{foxtoth}; and
Dvo\v{r}\'ak and Mohar's ingenious construction, 
for {each integer $k\ge 171$}, of
$k$-crossing-critical graphs of arbitrarily
large maximum degree~\cite{mohar}.

The earliest interesting, nontrivial construction of
$k$-crossing-critical graphs is due to \v{S}ir\'a\v{n}~\cite{siran}, who gave
examples of infinite families of $k$-crossing-critical graphs for
fixed values of $k$. These constructions involve graphs
with parallel edges. Shortly afterwards, Kochol~\cite{kochol} gave an infinite
family of $2$-crossing-critical, simple $3$-connected graphs.

In their influential paper on crossing-critical graphs, Richter and
Thomassen~\cite{richterthomassen} proved that $k$-crossing-critical
graphs have bounded crossing number. Richter and Thomassen also
investigated  regular simple crossing-critical
graphs. They used their aforementioned result to prove that for each
fixed $k$, there are only finitely many $k$-crossing-critical
$6$-regular simple graphs, and also constructed an
infinite family of $3$-crossing-critical, simple 
$4$-connected $4$-regular graphs.

We note that degree two vertices affect neither the crossing number
nor the crossing criticality of a graph. Also, the crossing number of
a disconnected graph is clearly the sum of the crossing numbers of its
connected components. Thus the interest in crossing-critical graphs
is focused on connected graphs with minimum degree at least $3$.

The construction of Richter and Thomassen  was generalized
in~\cite{salazar}, where it was shown that for every rational number
$q\in[4,6)$, there exists an integer $k_q$ such that there is an
infinite family of $k_q$-crossing-critical simple graphs with average
degree $q$.  Pinontoan and Richter~\cite{pinontoanrichter} extended
the range to every rational $q\in[3.5,6)$, and recently
Bokal~\cite{bokal} used his novel technique of zip products to
describe a construction that {yields} an infinite family for
 every rational $q\in(3,6)$.

What about $q = 3$ or $q \ge 6$?  Let $G$ and $H$ be simple $3$-regular
graphs.  Since $G$ has a subgraph isomorphic to a subdivision of $H$ 
if and only if $H$ is isomorphic to 
a minor of $G$, the Graph Minor Theorem~\cite{wagner} implies that 
for every integer $k\ge1$ there are only finitely many 
$k$-crossing-critical {3-regular} {simple}  {graphs}.
In fact, this does not need the full strength of the Graph Minor Theorem;
by Hlin\v{e}n\'y's result~\cite{petr} that $k$-crossing-critical graphs 
have bounded path-width all that is needed is the fact that
graphs of bounded path-width are well-quasi-ordered, which is a lot
easier that the general Graph Minor Theorem.
%
%
 On the other hand, it follows
easily from the techniques in~\cite{richterthomassen} that for each
fixed positive integer $k$ and rational $q>6$ there are only finitely
many $k$-crossing-critical simple graphs with average degree $q$.
Thus the only remaining open question is whether for some $k$ there exists
an infinite family of $k$-crossing-critical simple graphs of
average degree six. In this paper we answer this question in
the negative, as follows.

\begin{theorem}\label{thm:main}
For each fixed positive integer $k$, the collection of $k$-crossing-critical
simple graphs with average degree at least six is finite.
\end{theorem}

\ignore{\underline{DEFINITION for PLANE graphs.} 
The crucial new result behind the proof of Theorem~\ref{thm:main} is a statement
of independent interest. Let $G$ be a plane graph, and let 
$C_1,C_2,\ldots,C_t$\comm{Using $t$ instead of
  Robin's $k$} be pairwise distinct cycles in $G$. 
The sequence $C_1,C_2,\ldots,C_t$ is a {\em
  nest} if
for each $i=1,2,\ldots,t-1$, the cycle $C_{i+1}$ is contained in the
closed disk bounded by the cycle $C_i$.  We say that $t$ is the {\em
  size} of the nest.  If $X\subseteq V(G)$, $s:=|X|$\comm{Robin:
  used to be $s:=|X|\le 2$. Change ok?} and
$V(C_i)\cap V(C_j) = X$ for every two distinct indices
$i,j=1,2,\ldots,t$, then we say that $C_1,C_2,\ldots,C_t$ is an
$s$-nest.}

\noindent
In fact, we prove in {Theorem}~\ref{thm:main2} below that the conclusion holds for graphs of average degree
at least $6-c/n$, where $c$ is an absolute constant, and $n$ is the number
of vertices of the graph. 
{The assumption that $G$ be simple cannot be
omitted: as shown} in~\cite{richterthomassen}, 
{for each integer $p \ge
1$ there is an infinite family of
$4p$-regular $3p^2$-crossing-critical (nonsimple) graphs. Moreover, by
adding edges (some of them parallel) to the $4$-regular $3$-crossing-critical graphs
$H_m$} in~\cite{richterthomassen}, {it is possible to obtain an infinite
family of $6$-regular $12$-crossing-critical (nonsimple) graphs.}


The crucial new result behind the proof of Theorem~\ref{thm:main} is the
following theorem, which may be
 of independent interest. Let $\gamma$ be a planar drawing of a
graph $G$, and let $H$ be a subgraph of $G$. We say that $H$ is 
{\em  crossing-free in $\gamma$} if no edge of $H$ is crossed in $\gamma$
by another edge of $G$.
A sequence $C_1,C_2,\ldots,C_t$ of cycles in $G$
is a {\em nest in $\gamma$} if
the cycles are pairwise edge-disjoint, each of them is crossing-free in $\gamma$,
and for each $i=1,2,\ldots,t-1$ the cycle $\gamma(C_{i+1})$ is contained in
the closed disk bounded by  $\gamma(C_i)$.  We say that $t$ is the
{\em size} of the nest.  If $X\subseteq V(G)$, $s:=|X|$ and
$V(C_i)\cap V(C_j) = X$ for every two distinct indices
$i,j=1,2,\ldots,t$, then we say that $C_1,C_2,\ldots,C_t$ is an
$s$-{\em {nest}}.

\begin{theorem}
\label{label:nest}
For every integer $k$ there exists an integer $n$ such that every
planar triangulation on at least $n$ vertices has an $s$-nest of size
at least $k$ for some $s\in\{0,1,2\}$.
\end{theorem}

\noindent
To deduce Theorem~\ref{thm:main} from Theorem~\ref{label:nest} we prove
that a $k$-crossing-critical graph cannot have a large $s$-nest for
any $s\in\{0,1,2\}$.
For $s=2$ this was shown in~\cite{hlinenysalazar}, but we give
a shorter proof with a slightly better bound.

We formalize the notion of a planar drawing as follows. 
By a {\em polygonal arc} we mean a set $A\subseteq {\mathbb R}^2$
which is
the union of finitely many straight line segments and is homeomorphic
to the interval $[0,1]$. The images of $0$ and $1$ under the
homeomorphism are called the {\em ends} of $A$. A {\em polygon}
is a set $B\subseteq {\mathbb R}^2$ which is
the union of finitely many straight line segments and is homeomorphic
to the unit circle $\{(x,y)\in{\mathbb R}^2: x^2+y^2=1\}$.
Let $G$ be a graph. A {\em drawing} of $G$ is a mapping $\gamma$ 
with domain $V(G)\cup E(G)$ such that
\myitem{(i)}$\gamma(v)\in{\mathbb R}^2$ for every $v\in V(G)$,
\myitem{(ii)}$\gamma(v)\ne\gamma(v')$ for distinct $v,v'\in V(G)$,
\myitem{(iii)}for every non-loop edge $e\in E(G)$ with ends $u$ and $v$ there
   exists a polygonal arc $A\subseteq{\mathbb R}^2$ 
   with ends $\gamma(u)$ and $\gamma(v)$
   such that $\gamma(e)=A-\{u,v\}\subseteq{\mathbb R}^2-\gamma(V(G))$,
\myitem{(iv)}for every loop $e\in E(G)$ incident with $u\in V(G)$ there
   exists a polygon $P\subseteq{\mathbb R}^2$ containing $\gamma(u)$
   such that $\gamma(e)=P-\{u\}\subseteq{\mathbb R}^2-\gamma(V(G))$, and
\myitem{(v)}if $e,e'\in E(G)$ are distinct, then $\gamma(e)\cap \gamma(e')$
is finite.

\noindent
If $e,e'\in E(G)$ are distinct and $\gamma(e)\cap \gamma(e')\ne\emptyset$, 
then we say that $e$ and $e'$ {\em cross in $\gamma$} and that every
point of $\gamma(e)\cap \gamma(e')$ is a {\em crossing}.
(Thus a point where $\gamma(e)$ and $\gamma(e')$ ``touch" also counts
as a crossing.)
If $H$ is a subgraph of $G$, then
by $\gamma(H)$ we denote the image of $H$ under $\gamma$; that is,
the set of points in ${\mathbb R}^2$ that
either are equal to $\gamma(v)$ for some $v\in V(H)$ or belong to $\gamma(e)$
for some $e\in E(H)$.
A {\em plane graph} is a graph $G$ such that $V(G)\subseteq{\mathbb R}^2$,
every edge of $G$ is a subset of ${\mathbb R}^2$, and the identity mapping
$V(G)\cup E(G)\to V(G)\cup E(G)$ is a drawing of $G$ with no crossings.

We are restricting ourselves to piecewise linear drawings merely for convenience.
This restriction does not change the class of graphs that admit drawings
with a specified number of crossings, while piecewise linear drawings are
much easier to handle.

We prove Theorem~\ref{label:nest} in Section~\ref{structure}
and Theorem~\ref{thm:main} in Section~\ref{proofmain}.

\section{Finding a nest}
\label{structure}

A {\em tree decomposition} of a graph $G$ is a triple $(T,W,r)$
where $T$ is a tree, $r\in V(T)$ and $ W=(W_t~:~t \in V(T))$ is a collection of
subsets of $V(G)$ such that
\myitem{(T1)} $\bigcup_{t \in V(T)}W_t=V(G)$ and every edge of $G$ has both ends
in some $W_t$, and
\myitem{(T2)} if $t,t',t'' \in V(T)$ and $t'$ belongs to the unique path in $T$ connecting $t$ and $t''$, then $W_t \cap W_{t''} \subseteq W_{t'}$.

\noindent
The {\em width} of the tree-decomposition $(T,W,r)$ is the maximum of
$|W_t|-1$ over all $t\in V(T)$.
Now let $G$ be a plane graph.
We say that the tree-decomposition $(T,W,r)$ of $G$ is {\em standard} if
\myitem{(T3)} for every edge $e=tt'\in E(T)$ the set $W_t\cap W_{t'}$ is
the vertex-set of a cycle $C_e$ in $G$, and
\myitem{(T4)} if $e,e'\in E(T)$ are distinct, and $e$ lies on the unique
path from $r$ to $e'$, then $C_{e'}\ne C_e$
and $C_{e'}$ belongs to the closed disk bounded by $C_e$.

\noindent
The cycles $C_e$ will be called the {\em rings} of $(T,W,r)$.

We will need the following lemma.

\begin{lemma}
\mylabel{lem:stddec}
Let $k\ge1$ be an integer, and let $G$ be a triangulation of the plane.
Then $G$ has either a $0$-nest of size $k$, or a standard tree-decomposition
of width at most $12k-1$.
\end{lemma}

\proof
We may assume that $G$ has no $0$-nest of size $k$.
Let $(T,W,r)$ be a standard tree-decomposition of $G$ such that
\myitem{(a)} $T$ has at least one edge and maximum degree at most three,
\myitem{(b)} $|W_t|\le 12k$ if $t=r$ or if $t$ is not a leaf of $T$,
\myitem{(c)} each ring of $(W,T,r)$ has length at most $8k$,
\myitem{(d)} if $t\in V(T)-\{r\}$ and $t'$ is the unique neighbor of $t$
in the path in $T$ from $t$ to $r$, then $W_t$ consists precisely of
the vertices of $G$ drawn in the closed disk bounded by $C_{tt'}$, and

\noindent
subject to (a)--(d)
\myitem{(e)} $T$ is maximal.

\noindent 
Such a choice is possible, because of the following construction.
Let $T$ be a tree with vertex-set $\{t_1,t_2\}$, let $C$ be the
triangle bounding the outer face of $G$, let $W_{t_1}=V(C)$,
and let $W_{t_2}=V(G)$. Then $(T,W,t_1)$ satisfies (a)--(d).

So let $(T,W,r)$ satisfy (a)--(e).
We claim that $(T,W,r)$ has width at most $12k-1$.
To prove that suppose to the contrary that $|W_{t_0}|>12k$ for
some $t_0\in V(T)$.
Then by (b) $t_0\ne r$ and $t_0$ is a leaf of $T$.
Let $t_1$ be the unique neighbor of $t_0$ in $T$, and let
$C$ denote the ring $C_{t_0t_1}$.
Then $|V(C)|\le 8k$ by (c).
Let $\Delta$ denote the closed disk bounded by $C$, and let $H$
be the near-triangulation consisting of all vertices and edges of
$G$ drawn in $\Delta$.
By (d) we have $V(H)=W_{t_0}$.
For $u,v\in V(C)$, let $c(u,v)$ (respectively, $d(u,v)$) be the number of
edges in the shortest path of $C$ (respectively, $H$) between $u$ and $v$.

\myclaim{1}{$c(u,v)=d(u,v)$ for all $u,v\in V(C)$.}

To prove (1) we certainly have
$d(u,v)\le c(u,v)$ since $C$ is a subgraph of $H$.
If possible, choose a pair $u,v\in V(C)$ with $d(u,v)$ minimum such that
$d(u,v)<c(u,v)$. Let $P$ be a path of $H$ between $u$ and $v$, with
$d(u,v)$ edges.  Suppose that some internal vertex $w$ of $P$ belongs
to $V(C)$.  Then
$$d(u,w)+d(w,v)=d(u,v)<c(u,v)\le c(u,w)+c(w,v)$$
and so either $d(u,w)<c(u,w)$ or $d(w,v)<c(w,v)$, in either case contrary
to the choice of $u,v$. Thus there is no such $w$. Let
$C,C_1,C_2$ be the three cycles of $C\cup P$, let
$\Delta, \Delta_1,\Delta_2$ be the closed disks they bound, and for $i=1,2$
let $H_i$ be the subgraph of $H$ consisting of all vertices and edges
drawn in $\Delta_i$.
Then $C_1$ and $C_2$ have length at most $8k$.
Let $T'$ be the tree obtained from $T$ by adding two vertices $r_1,r_2$,
both joined to $t_0$.
For $t\in V(T)-\{t_0\}$ let $W'_{t}=W_t$, let $W'_{t_0}=V(C\cup P)$,
let $W'_{r_i}=V(H_i)$, and let $W'=(W'_{t}: t\in V(T'))$.
Then $(T',W',r)$ satisfies (a)-(d), contrary to (e).
This proves (1).

\myclaim{2}{$C$ has length exactly $8k$.}

To prove (2) suppose for a contradiction that $C$ has length at most $8k-1$.
Let $uv$ be an edge of $C$, and let $w$ be the third vertex of the face
incident with $uv$ and contained in the disk bounded by $C$.
Then $w\not\in V(C)$ by (1) and the fact that $|W_{t_0}|>12k$.
Let $T'$ be obtained from $T$ by adding a new vertex $r_0$ joined to $t_0$,
for $t\in V(T)-\{t_0\}$ let $W_t'=W_t$, let $W'_{t_0}=V(C)\cup\{w\}$,
let $W'_{r_0}=W_{t_0}$, and let $W'=(W_t':t\in V(T'))$.
Then $(T',W',r)$ is a standard tree-decomposition satisfying (a)--(d),
contrary to (e). This proves (2).

Now let $v_1,v_2,\ldots,v_{8k}$ be the vertices of $C$ in order.
By (1) and~\cite[Theorem (3.6)]{RobSeyGM6} there exist $2k$ disjoint paths
from $\{v_1,v_2,\ldots,v_{2k}\}$ to $\{v_{4k+1},v_{4k+2},\ldots,v_{6k}\}$,
and $2k$ disjoint paths
from $\{v_{2k+1},v_{2k+2},\ldots,v_{4k}\}$ to $\{v_{6k+1},v_{6k+2},\ldots,v_{8k}\}$.
Using those sets of paths it is easy to construct a $0$-nest in $G$ of size $k$.
In fact, using the argument of~\cite[Theorem (4.1)]{RobSeyGM3} it can be shown
that $G$ has a $2k\times 2k$ grid minor, and hence a $0$-nest of size $k$.~\qed
\bigskip

\noindent
{\bf Proof of Theorem~\ref{label:nest}.}
Let $k\ge1$ be a given integer, 
let $h$ be an integer such that for every coloring
of the edges of the complete graph on $h$ vertices using at most $12k$
colors, there is a monochromatic clique of size $24k^2$,
and let $n=36k\cdot2^{h+1}$.
The integer $h$ exists by Ramsey's theorem.
We claim that $n$ satisfies the conclusion of the theorem.
To prove the claim let $G$ be a triangulation of the plane on at least
$n$ vertices.
By Lemma~\ref{lem:stddec} we may assume that $G$ has a standard
tree-decomposition $(T,W,r)$ of width at most $12k$.
It follows that $T$ has at least $n/(12k)$ vertices.
Thus $|V(T)|> 3\cdot 2^{h+1}-2$, and hence $T$ has a path on
$h+1$ vertices starting in $r$.
Let $t_0=r,t_1,\ldots,t_h$ be the vertices of one such path,
and for $i=1,2,\ldots,h$ let $C_i$ denote the ring $C_{t_{i-1}t_i}$.
Then by (T3) and (T4) $C_1,C_2,\ldots,C_h$  is a sequence  of distinct cycles
such that for indices $i,j$ with $1\le i\le j\le h$ the cycle $C_j$
belongs to the closed disk bounded by $C_i$.
We shall refer to the latter condition as the nesting property.
Let $K$ be a complete graph with vertex-set $\{1,2,\ldots,h\}$.
We color the edges of $K$ by saying that the edge $ij$ is colored
using $|V(C_i)\cap V(C_j)|$.
By the choice of $n$ there exist a subsequence $D_1,D_2,\ldots,D_{24k^2}$
of $C_1,C_2,\ldots,C_h$ and an integer $t\in\{0,1,\ldots,12k-1\}$
such that $|V(D_i)\cap V(D_j)|=t$ for every pair of distinct integers
$i,j\in\{1,2,\ldots,24k^2\}$.
Since the sequence $D_1,D_2,\ldots,D_{24k^2}$ satisfies the nesting property, 
we deduce that there exists a set $X$
such that $V(D_i)\cap V(D_j)=X$ for every pair of distinct integers
$i,j\in\{1,2,\ldots,24k^2\}$.
If $|X|\le1$, then the sequence $D_1,D_2,\ldots,D_{k}$ satisfies the
conclusion of the theorem.
We may therefore assume that $|X|\ge2$.
Let the elements of $X$ be numbered $x_1,x_2,\ldots,x_t=x_0$ in such a way
that they appear on $D_1$ in the order listed. 
It follows that they appear on each cycle $D_j$ in the order listed.
Now for $i=1,2,\ldots,t$ and $j=1,2,\ldots,24k^2$ let $P_{ij}$ be 
the subpath of $D_j$ with ends $x_{i-1}$ and $x_i$ that is disjoint from 
$X-\{x_{i-1},x_i\}$ (if $|X|=2$ we number the two subpaths of $D_j$ arbitrarily).
Since the cycles $D_j$ are pairwise distinct and $t\le 12k-1$, 
we deduce that there exists an index $i\in\{1,2,\ldots,t\}$
such that the path $P_{ij}$ has at least one internal vertex for
at least $2k$ distinct integers $j\in\{1,2,\ldots,24k^2\}$.
Let us fix this index $i$, and let $Q_1,Q_2,\ldots,Q_{2k}$
be a subsequence of $P_{i1},P_{i2},\ldots$ such that each $Q_j$ has
at least one internal vertex.
It follows that the paths $Q_1,Q_2,\ldots,Q_{2k}$ are internally disjoint
and pairwise distinct.
Thus $Q_1\cup Q_{2k}, Q_2\cup Q_{2k-1},\ldots,Q_k\cup Q_{k+1}$ is
a $2$-nest in $G$ of {size} $k$, as desired.~\qed

\section{Using a nest}
\label{proofmain}

To prove Theorem~\ref{thm:main} we need several lemmas,
but first we need a couple of definitions.
We say that an $s$-nest
$C_1,C_2,\ldots,C_t$ in a drawing $\gamma$ of a graph $G$ is {\em clean} if
every crossing in $\gamma$ belongs either to the open disk bounded by
$\gamma(C_t)$, or to the complement of the closed disk bounded by $\gamma(C_1)$.
We say that a drawing $\gamma$ of a graph $G$ is {\em generic} if
it satisfies (i)-(v) and
\myitem{(vi)}every point $x\in{\mathbb R}^2$ belongs to $\gamma(e)$ for
at most two edges $e\in E(G)$, and
\myitem{(vii)}if $\gamma(e)\cap \gamma(e')\ne\emptyset$ for distinct
edges $e,e'\in E(G)$, then $e$ and $e'$ are not adjacent.

\def\ka{{\ell}}

\begin{lemma}
\label{nestexists}
For every three integers {$\ka$}, $r,t\ge0$  there exists 
an integer $n_0$ such that for every simple graph $G$ on $n\ge n_0$ vertices
of average degree at least $6-r/n$ and 
every generic drawing $\gamma$ of $G$ with at most {$\ka$} crossings
there exists an $s$-nest in $\gamma$ of size $t$ for some $s\in\{0,1,2\}$.
\end{lemma}

\proof
Let $\ka,t,r$ be given, and let $n_0$ be an integer such that
Theorem~\ref{label:nest} holds when {$k$} is replaced by 
$t':=t+2\ka+r-6$ and $n$ is replaced by $n_0$.
We will prove that $n_0$ satisfies the conclusion of the theorem.
To that end let $G$ be a simple graph on $n\ge n_0$ vertices
of average degree at least $6-r/n$ and
let $\gamma$ be a generic drawing of $G$ with at most $\ka$ crossings.
We will prove that $\gamma$ has a desired $s$-nest.
%
Let $G'$ denote the plane graph obtained from $\gamma$ by converting
each crossing into a vertex.
Let $V_4$ be the set of these new vertices. Then $|V_4|\le \ka$.
By (vi) each vertex in $V_4$ has degree four in $G'$,
and since $G$ is simple it follows from (vii) that $G'$ is simple.
Let $\deg(v)$ denote the degree of $v$ in $G'$, let $\cal F$ denote
the set of faces of $G'$, and for $f\in{\cal F}$ let $|f|$ denote
the length of the boundary of $f$; that is, the sum of the lengths of the
walks forming the boundary of $f$.
By Euler's formula we have
$$
  \sum_{v\in V(G')}(6-\deg(v)) + \sum_{f\in\cal F} 2(3-|f|) = 12.
$$
But $\sum_{v\in V(G)-V_4}(6-\deg(v))\le r$ by hypothesis, and so
$$
   \sum_{f\in\cal F} (|f|-3)\le {1\over 2}\sum_{v\in V_4}(6-\deg(v))-6+r=
   |V_4|-6+r\le \ka+r-6,
$$
because every vertex in $V_4$ has degree four in $G'$.
Thus $G'$ has at most $\ka+r-6$ non-triangular faces, each of size at
most $\ka+r-3$.

Let $G''$ be the triangulation obtained from $G'$ by adding a vertex into each
non-triangular face and joining it to each vertex on the boundary of that face.
Thus every added vertex has degree in $G''$ at most $\ka+r-3$. 
By Theorem~\ref{label:nest} the triangulation $G''$ has an $s$-nest
$C_1,C_2,\ldots,C_{t'}$ of {size} $t'$ for some $s\in\{0,1,2\}$.
Let $X$ be the set of vertices every two distinct cycle $C_i$ and $C_j$
have in common.
Then every vertex of $X$ has degree at least $2t'$, and hence belongs
to $V(G)$, because every vertex of $V(G'')-V(G)$ has degree four
or at most $\ka+r-3$.
Thus at most $\ka+(\ka+r-6)=2\ka+r-6$ cycles $C_i$ contain a vertex not
in $G$, and by removing all those cycles we obtain a desired $s$-nest in
$\gamma$.~\qed

\begin{lemma}
\label{lem:clean}
Let {$k\ge 0$} and $t\ge1$ be  integers,  let $s\in\{0,1,2\}$, let $G$ be a graph,
and let $\gamma$ be a drawing of  $G$ with at most {$k$} crossings and
an $s$-nest of size {$(k+1)$} $(t-1)+1$. 
Then $\gamma$ has a clean $s$-nest of size $t$.
\end{lemma}

\proof
Let $C_1,C_2,\ldots,C_{(k+1)(t-1)+1}$ be an $s$-nest in $G$. 
For $i=1,2,\ldots,(k+1)(t-1)$ let $\Omega_i$ denote the subset of ${\mathbb R}^2$
obtained from the closed disk
bounded by $\gamma(C_i)$ by removing the open disk bounded by $\gamma(C_{i+1})$.
Since there are at most $k$ crossings in $\gamma$, it follows that
$\Omega_i$ includes no crossing of $\gamma$ for $t-1$ consecutive
integers  {in} $\{1,2,\ldots,(k+1)(t-1)\}$, say $i,i+1,\ldots,i+t-2$.
Then $C_i,C_{i+1},\ldots,C_{i+t-1}$ is a clean $s$-nest of {size} $t$,
as desired.~\qed


\begin{lemma}
\label{lem:cr}
Let $k\ge1$ be an integer, let $s\in\{0,1\}$ and let $\gamma$ be a
drawing of a graph $G$ with a clean $s$-nest of size $4k+1$.
If $\hbox{\rm cr}(G-e)< k$ for all $e\in E(G)$, then $\hbox{\rm cr}(G)< k$.
\end{lemma}

\proof
If $C$ is a cycle of $G$ that is crossing-free in $\gamma$, then we
denote by $\Delta(C)$ the disk bounded by $\gamma(C)$.
Let $D_1,D_2,\ldots,D_{4k+1}$ be a clean $s$-nest in $\gamma$ of size $4k+1$.
We may assume that the $s$-nest is chosen so that
\myitem{(1)}for $i=2,3,\ldots,2k$, if $D$ is a cycle in $G$ such that
  $D_{i-1},D,D_{i+1}$ is an $s$-nest in $\gamma$ and
  $\Delta(D_i)\subseteq\Delta(D)$, then $D_i=D$, and
\myitem{(2)}for $i=2k+2,2k+3,\ldots,4k$, if $D$ is a cycle in $G$ such that
  $D_{i-1},D,D_{i+1}$ is an $s$-nest in $\gamma$ and
  $\Delta(D)\subseteq\Delta(D_i)$, then $D_i=D$.

\noindent
Let $e\in E(D_{2k+1})$.
By hypothesis there exists a drawing $\gamma'$ of $G-e$ with at most
$k-1$ crossings.
Thus at most $2k-2$ cycles among $D_1,D_2,\ldots,D_{2k}$ include an
edge that is crossed by another edge in $\gamma'$, and
similarly for $D_{2k+2},D_{2k+3},\ldots,D_{4k+1}$.
Hence there exist indices $i_2\in\{2,3,\ldots,2k\}$ and
$i_4\in\{2k+2,2k+3,\ldots,4k\}$ such that $D_{i_2}$ and $D_{i_4}$
are  crossing-free in $\gamma'$.
Let $C_1:=D_{i_2-1}$, $C_2:=D_{i_2}$, $C_3:=D_{2k+1}$,
 $C_4:=D_{i_4}$, and $C_5:=D_{i_4+1}$.
Then $C_1,C_2,C_3,C_4,C_5$ is a clean $s$-nest in $\gamma$.
Let $H:=C_2\cup C_4$, let $B_1$ be the $H$-bridge containing $C_1$,
and let $B_5$ be the $H$-bridge containing $C_5$.
Since $H$ is crossing-free in $\gamma$ we see that $B_1\ne B_5$.
Let $\Omega$ be the face of $\gamma(H)$ that is incident with edges of
both $C_2$ and $C_4$.
Thus if $s=0$, then $\Omega$ is an annulus, and if $s=1$ it is
a ``pinched annulus".
An $H$-bridge $B$ of $G$ will be called {\em interior} if
$\gamma(B)$ is a subset of the closure of $\Omega$ and it will be
called {\em exterior} otherwise.
We need the following claim.

\myclaim{3}{If $B$ is an exterior $H$-bridge of $G$ with at least two
   attachments, then either $B=B_1$, or $B=B_5$.}

To prove (3) let $B$ be an exterior $H$-bridge of $G$ with at least two
attachments.
Since $C_2$ and $C_4$ are crossing-free in $\gamma$ it follows that
either all attachments of $B$ belong to $C_2$, or they all belong to $C_4$.
From the symmetry we may assume the former.
We may also assume that $B\ne B_1$, for otherwise the claim holds.
Thus $B$ is disjoint from $C_1$.
Since $B$ has at least two attachments it includes a path $P$ with both
ends in $C_2$ and otherwise disjoint from it.
Since $C_1,C_2,C_3,C_4,C_5$ is a clean $s$-nest, no edge of $P$ is
crossed by another edge in $\gamma$.
Thus $C_2\cup P$ includes a cycle $D$ disjoint from $C_1=D_{i_2-1}$
with $\Delta(C_2)\subseteq\Delta(D)$, contrary to (1).
This proves (3).
\medskip

We may assume, by composing $\gamma'$ with a homeomorphism of the plane,
that $\gamma(H)=\gamma'(H)$.
We now define a new drawing $\delta$ of $G$ as follows.
For every vertex and edge $x$  that belongs to $H$ or to an interior
$H$-bridge of $G$ we define $\delta(x)=\gamma(x)$.
If $\gamma'(B_5)\subseteq\Delta(C_4)$, then for every
$x\in V(B_5)\cup E(B_5)$ we define $\delta(x)=\gamma'(x)$;
otherwise we use circular inversion to redraw $\gamma'(B_5)$ in
$\Delta(C_4)$ and use the inversion of $\gamma'(B_5)$ to define $\delta(B_5)$.
We define $\delta(B_1)$ analogously.
Finally, for an $H$-bridge $B$ with at most one attachment we define
$\delta(B)$ by scaling $\gamma'(B)$ suitably so that it does not
intersect any other $H$-bridge of $G$.
Thus every crossing of $\delta$ is also a crossing of $\gamma'$,
and hence $\delta$ has at most $k-1$ crossings, as desired.~\qed
\bigskip

We also need a version of Lemma~\ref{lem:cr} for $2$-nests.
Such a lemma follows from~\cite[Theorem~1.3]{hlinenysalazar},
but we give a proof from first
principles, because we have already done a lot of the needed work
in the previous lemma.
Furthermore, Lemmas~\ref{lem:clean} and~\ref{lem:cr2} imply
a small numerical improvement to~\cite[Theorem~1.3]{hlinenysalazar}.

\begin{lemma}
\label{lem:cr2}
Let $k\ge1$ be an integer, and let $\gamma$ be a
drawing of a graph $G$ with a clean $2$-nest of size $4k+1$. 
{If $\hbox{\rm cr}(G-e)< k$ for all $e\in E(G)$, then $\hbox{\rm cr}(G)< k$.}
\end{lemma}

\proof
We proceed similarly as in Lemma~\ref{lem:cr}.
Let $\Delta(C)$ be as in Lemma~\ref{lem:cr}.
We select our $2$-nest so that it satisfies (1) and (2) from
Lemma~\ref{lem:cr}.
We pick $e\in E(D_{2k+1})$, but we now require that $e$ be incident
with a vertex in $X$.
{Let $\gamma'$ be a drawing of $G-e$ with at most $k-1$ crossings.}
We choose a clean nest $C_1,C_2,C_3,C_4,C_5$ in the same way as before,
with one caveat: the index $i_4$ can be chosen so that there exists
an index $i_5$ such that $i_4<i_5\le 4k+1$ and $D_{i_5}$ is crossing-free
in $\gamma'$.
We put $C_6:=D_{i_5}$.
Now $\gamma(H)$ has four faces, and we define $\Omega$ to be the one
containing $\gamma(e)$.
For $i=1,2,\ldots,5$ let $P_i$ be a subpath of $C_i$ with ends in $X$
chosen as follows: $P_2$ and $P_4$ are defined by saying that
$P_2\cup P_4$ is the boundary of $\Omega$, $P_3$ is defined by
$\gamma(P_3)\subseteq\Omega$, and $P_1$ and $P_5$ are defined by saying that
$\Delta(P_1\cup P_2)$ and $\Delta(P_4\cup P_5)$ are disjoint from $\Omega$.
We define $B_1$ as the $H$-bridge of $G$ containing $P_1$
and $B_5$ as the $H$-bridge containing $P_5$.
It is now possible that $B_1=B_5$.

Let us say that an $H$-bridge is {\em singular} if its set of attachments
is $X$. The following is an analogue of claim (3) from Lemma~\ref{lem:cr}.
The proof follows the same lines, and so we omit it.

\myclaim{3}{If $B$ is an exterior $H$-bridge of $G$ with at least two
   attachments, then either $B=B_1$, or $B=B_5$, or $B$ is singular.}

Again, we may assume that $\gamma(H)=\gamma'(H)$.
If $B_1\ne B_5$, then the argument proceeds in the same way as in
Lemma~\ref{lem:cr}. All singular $H$-bridges can be drawn outside $\Omega$
so that they will be disjoint from each other and from
all other $H$-bridges.

Thus we may assume that $B_1=B_5$.
Now we may assume, by applying a circular inversion with respect to $H$
if necessary,
that $\gamma'(B_1)$ lies in the complement of $\Omega$.
Finally, we may assume that, subject to the conditions already imposed
on $\gamma'$,

\myclaim{4}{the number of $H$-bridges $B$ with $\gamma'(B)$ contained
   in the closure of $\Omega$ is minimum.}

Let $d_1$ be the maximum number of edge-disjoint paths in $G-X$
from $P_2$ to $P_4$ that are contained in interior $H$-bridges,
and let $d_2$ be the maximum number of edge-disjoint such paths
contained in exterior $H$-bridges.
Let $d_1',d_2'$ be defined analogously, but with respect to the
drawing $\gamma'$.

\myclaim{5}{$d_1=d_1'$ and $d_2=d_2'$}

To prove (5) we notice that a path in $G-X$ from $P_2$ to $P_4$ belongs
to a $H$-bridge that has an attachement in both $P_2-X$ and $P_4-X$.
Let us call such an $H$-bridge {\em global}.
Then $B_1$ is the only global exterior bridge by (3),
and $\gamma'(B_1)\cap \Omega=\emptyset$, because we chose $\gamma'$ that way.
Conversely, if $B$ is a global interior bridge, then
$\gamma'(B-e)$ lies in the closure of $\Omega$, for otherwise $\gamma'(B-e)$
would have to intersect $\gamma'(C_6)$ (because $C_6$ is
crossing-free in $\gamma'$), and hence $B=B_1$, contrary
to the fact that $B$ is interior.
Since $e$ is incident with $X$, claim (5) follows.

\myclaim{6}{If there exists an exterior singular $H$-bridge, then $d_1\ge d_2$.}

To prove (6) let $B$ be an exterior singular bridge.
Let $d$ be the maximum number of edge-disjoint paths in $B$ that
join the vertices of $X$.
Then $d>0$, because $B$ is singular.
In $\gamma$ the bridges $B$ and $B_1$ cross at least $dd_2$ times
by definition of $d$ and $d_2$.
By Menger's theorem there is a set of edges $F\subseteq E(B)$
of size $d$
such that $B-F$ includes no path joining the two vertices in $X$.
Let $J$ be the union of $H$ and all interior $H$-bridges.
Similarly, there is a set of edges $F_1\subseteq E(J)$ of size
$d_1$ such that $J-F_1$ has no edge from $P_2$ to $P_4$.
Now $\gamma$ can be changed by redrawing $B$ inside $\Omega$.
This can be done in such a way that every edge of $F$ crosses every
edge of $F_1$, and these are the only crossings of an edge of $B$
with an edge not in $B$.
Thus the redrawing of $B$ removes the at least $dd_2$ crossings
of $B$ and $B_1$ and introduces exactly $dd_1$ new crossings.
Since $G$ has crossing number exactly $k$, the new drawing has at least
$k$ crossings, and hence $d_1\ge d_2$, as desired.
This proves (6).

\myclaim{7}{Every  exterior $H$-bridge $B$ satisfies
    $\gamma'(B)\cap\Omega=\emptyset$.}

To prove (7) let $B$ be an exterior $H$-bridge.
If $B$ has only one attachment, then it clearly satisfies the conclusion
of (7), for otherwise it can be redrawn outside of $\Omega$,
contrary to (4).
Thus $B$ has at least two attachments, and hence $B=B_1$, or $B$
is singular by (3). If $B=B_1$, then the conclusion of (7) holds, and so
we may assume that $B$ is singular.
By (6) $d_1\ge d_2$, and hence $d_1'\ge d_2'$ by (5).
If $\gamma'(B)\cap\Omega\ne\emptyset$, then using the argument of (6)
we can change the drawing $\gamma'$ by drawing $B$ in the complement
of $\Omega$, contrary to (4).
This proves (7).
\medskip

Let $\delta$ be a drawing of $G$ defined to coincide with $\gamma$
on $H$ and all interior $H$-bridges, and to coincide with $\gamma'$ for
all other $H$-bridges.
By (7) this is a well-defined drawing of $G$, and every crossing of
$\delta$ is a crossing of $\gamma'$, because there are no crossings
in $\Omega$.
Thus $\delta$ has at most $k-1$ crossings, 
{as desired}.~\qed
\bigskip

We are now ready to prove Theorem~\ref{thm:main}, which we restate
in a slightly stronger form.

\begin{theorem}
\label{thm:main2}
For all integers {$k \ge 1$}, $r\ge0$ {there is an integer $n_0:=n_0(k,r)$ such
that if $G$ is a $k$-crossing-critical simple graph on $n$ vertices
with average degree at least $6-r/n$, then $n< n_0$.} 
\end{theorem}

\proof Let {$k\ge 1$}, $ r\ge0$ be integers, let 
{$\ka := 2.5k+16$}, let 
{$t:=4k(\ell+1)+1$}, let $n_0$
be an integer such that Lemma~\ref{nestexists} holds, and
let $G$ be a $k$-crossing-critical
simple graph on $n$ vertices with average degree at least $6-r/n$.
We claim that $n<n_0$.
To prove the claim suppose to the contrary that $n\ge n_0$, and let
$\gamma$ be a drawing of $G$ with at most {$\ka$ crossings; such a
drawing exists by}~\cite[Theorem 3]{richterthomassen}.
By a standard and well-known argument we may assume that $\gamma$ is generic.
By Lemma~\ref{nestexists} there is an integer $s\in\{0,1,2\}$ and
an $s$-nest in $\gamma$ of size $t$,
and by Lemma~\ref{lem:clean} there is a clean $s$-nest in $\gamma$ of size $4k+1$.
That contradicts Lemma~\ref{lem:cr} if $s\in\{0,1\}$, or Lemma~\ref{lem:cr2}
if $s=2$.
Thus $n<n_0$, {as desired}.~\qed

\baselineskip 11pt
\vfill
\noindent
This material is based upon work supported by the National Science Foundation.
Any opinions, findings, and conclusions or
recommendations expressed in this material are those of the authors and do
not necessarily reflect the views of the National Science Foundation.
\eject

\end{document}